\documentclass[12pt]{article} 
\usepackage[reqno]{amsmath} 
\usepackage{amssymb, amsthm}
\usepackage{enumerate} 
\usepackage{graphicx}
\usepackage{microtype}

\newtheoremstyle{plainsl}%
	{\topsep}
	{\topsep}
	{\slshape} 
	{}
	{\normalfont\bfseries}
	{.}
	{ }
	{}

\theoremstyle{plainsl}
\newtheorem{theorem}{Theorem}[section]
\newtheorem{lemma}[theorem]{Lemma}
\newtheorem{corollary}[theorem]{Corollary}

\renewcommand\proof{\noindent\textsl{Proof. }}
\newcommand\sqr[2]{{\vbox{\hrule height.#2pt
    \hbox{\vrule width.#2pt height#1pt \kern#1pt
        \vrule width.#2pt}\hrule height.#2pt}}}

\renewcommand\qed{%
	\ifmmode\eqno\sqr53
	\else\nolinebreak\ \hfill\sqr53\medbreak\fi}

\numberwithin{equation}{section}

\newcommand\cP{{\mathcal P}}

\newcommand\cT{{\mathcal T}}

\title{Partially ordering the class of invertible trees}
\author{Krystal Guo\thanks{D\'{e}partment de Math\'{e}matique, Universit\'{e} libre de Bruxelles, Brussels, Belgium. \texttt{guo.krystal@gmail.com} Part of this work was done when the author was a post-doctoral fellow at University of Waterloo.}}

\begin{document}
\maketitle

\begin{abstract}

A tree $T$ is invertible if and only if $T$ has a perfect matching.  In \cite{Go85}, Godsil considers an invertible tree $T$ and finds that the matrix $A(T)^{-1}$ has entries in $\{0,\pm1\}$ and is the signed adjacency matrix of a graph which contains $T$. In this paper, we give a new proof of this theorem, which gives rise to a partial ordering relation on the class of all invertible trees on $2n$ vertices. In particular, we show that given an invertible tree $T$ whose inverse graph has strictly more edges, we can remove an edge from $T$ and add another edge to obtain an invertible tree $\widetilde{T}$ whose median eigenvalue is strictly greater. This extends naturally to a partial ordering. We find the maximal and minimal elements of this poset and explore the implications about the median eigenvalues of invertible trees.

\end{abstract}

\section{Introduction}

A tree is a connected graph with no cycles. It is clear from the expansion of the determinant that the adjacency matrix of a tree $T$ is invertible if and only if $T$ has a perfect matching. Such a tree is said to be \textsl{invertible}. In \cite{Go85}, Godsil considers an invertible tree $T$ and finds that the matrix $A(T)^{-1}$ has entries in $\{0,\pm1\}$ and is the signed adjacency matrix of a graph which contains $T$. We will refer the underlying graph of $A(T)^{-1}$ as the \textsl{inverse graph} of an invertible tree. 

For a graph $X$, let $\lambda_i(X)$ denote the $i$th largest eigenvalue of $X$. For a graph $X$ on $m$ vertices, the median eigenvalues of $X$ are $\lambda_{\lceil m/2 \rceil}$ and $\lambda_{\lfloor m/2 \rfloor}$. Median eigenvalue of bipartite graphs have been studied in \cite{Moh2015a,Moh2015b}. Median eigenvalues have applications in mathematical chemistry as they are related to the HOMO-LUMO separation, see for example \cite{FowPis2010}.

Since the eigenvalues of a bipartite graph are symmetric about $0$, the median eigenvalues of a bipartite graph $X$ on $2m$ vertices are such that $\lambda_m = -\lambda_{m+1}$. We will consider only the median eigenvalue for bipartite graph on an even number of vertices, since the median eigenvalue of any bipartite graph on an odd number of vertices is $0$ and it will suffice to consider only the positive median eigenvalue. For an invertible tree, we see that $\lambda$ is an eigenvalue of $T$ if and only if $ 1/\lambda$ is an eigenvalue of $A(T)^{-1}$. In \cite{Go85}, Godsil characterises the invertible trees on $2n$ vertices which attain the smallest $\lambda_n$; he proved for any tree $T$ on $2n$ vertices
\[\lambda_n(T) \geq \lambda_n(P_n) \]
and that equality is attained if and only if $T \cong P_n$. 

Subsequently classifications of trees attaining the maximum median eigenvalue has been found in \cite{BarNeuPat2006}, which also characterizes the trees which are isomorphic to their inverse graphs. A generalization for bipartite graphs is studied in \cite{SimCao89}. A characterization of all bipartite graphs with a
unique perfect matching whose adjacency matrices have inverses diagonally similar to non-negative matrices is given in \cite{YanYe17}. Graph inverses are an interesting topic and have  also studied in \cite{Tif11, McLMcN14, PanPat2016, YeYangManKle2017}.

In this paper, we give a new proof of Godsil theorem on the inverse of trees, using newer techniques. In the process, we define a partial ordering relation on the class of all invertible trees on $2n$ vertices. In particular, given an invertible tree $T$ whose inverse graph has strictly more edges than $T$, we can obtain a non-isomorphic invertible tree $\widetilde{T}$ such that $\lambda_n(T) < \lambda_n(\widetilde{T})$. This extends naturally to a partial ordering. We show that the maximal elements of this poset are exactly the trees obtained as the rooted product of a tree on $n$ vertices with $K_2$ and  the minimal elements of this poset are the path graphs.  

\section{Preliminaries}
\label{sec:pre} 

In this section, we will state results in the literature which will allow us to consider the inverse of a tree. We will denote the adjacency matrix of a graph $X$ by $A(X)$, or $A$ when the context is clear. If $X$ is a graph such that $A(X)$ is invertible, we say $X$ itself is \textsl{invertible} for brevity. We will refer to standard texts such as \cite{G93} and \cite{CDS95} for further background. 

The following theorem can be derived from Lemma 2.1 in Chapter 4 of \cite{G93}. 

\begin{theorem} If $a$ and $b$ are distinct vertices in $X$, then
\begin{equation}\label{eq:ab}
(tI -A)^{-1}_{a.b} = \frac{\phi_{a,b}(X,t)}{\phi(X,t)} = \frac{\sum_{P \in \cP_{a,b}}\phi(X\setminus P ,t)}{\phi(X,t)}
\end{equation}
where $P_{a,b}$ denotes the set of all $ab$-paths in $X$ and $X \setminus P$ is the graph obtained from $X$ by deleting all vertices of $P$.
\end{theorem}

For the diagonal entries we will use:

\begin{theorem} If $a$ is a vertex of  $X$, then
\begin{equation}\label{eq:aa} 
(tI -A)^{-1}_{a.a} = \frac{\phi(X\setminus a, t)}{\phi(X,t)}.
\end{equation}
\end{theorem}

Let $m(X)$ denote the number of perfect matching of a tree $T$. If $A(T)$ is invertible, then $[t^0] \phi(T,t) \neq 0$ and $T$ has a perfect matching and necessarily has an even number of vertices.  Let $T$ be a tree on $2n$ vertices. Since every tree has at most $1$ perfect matching, we obtain that 
\[
[t^0] \phi(T,t) = \begin{cases} (-1)^{n}, &\text{ if } m(T) =1; \\
0, &\text{ otherwise}. \end{cases}
\] 
In a tree, there is a unique path between any two distinct vertices $a$ and $b$; we will denote this path by $P_{a,b}$. If $M$ is a matching in $X$, we will say a path $P$ is \textsl{$M$-alternating path} if $P$ has $2k-1$ edges and contains $k$ edges of $M$.

\begin{lemma}\label{lem:alt} Let $T$ be an invertible tree and let $M$ be the unique perfect matching.  If $a,b$ be distinct vertices of $T$, then $T\setminus P_{a,b}$ has a perfect matching if and only if $P_{a,b}$ is a $M$-alternating path.  \end{lemma}

\proof If $ab$ is an edge of $X$, then $X\setminus\{a,b\}$ has a perfect matching if and only if $ab \in M$. 

Suppose $X\setminus P_{a,b}$ has a perfect matching $M'$. Consider the graph $Y = (V(X), M \Delta M')$ where $\Delta$ denotes symmetric difference. In $Y$, the vertices of $P_{a,b}$ have degree $1$. Every other vertex is incident to an edge of both $M$ and $M'$ and so their degree in $Y$ is either $2$ or $0$. Thus, $P_{a,b}$ is a path in $X$ which alternates between edges of $M$ and edges not in $M$ and (since $a,b$ are covered by $M$) has odd length, say $2k-1$ and has $1$ more edge in $M$ than not in $M$. 

For the converse, let $P_{a,b}$ be an $M$-alternating path of length $2k-1$ in $X$ with end vertices $a$ and $b$, which contains $k$ edges of $M$. If we take the edges of $M \setminus (P_{a,b}\cap M)$, we obtain a matching of $X\setminus P_{a,b}$. \qed 

\section{Inverses of trees}\label{sec:inv-tree} 

We will use the material in the previous section to study inverses of trees. We use different tools from those used in the original proof of Godsil's theorem about the inverses of trees and will give an alternate proof for his theorem. This will also lay the groundwork for the partial order on the class of invertible trees in Section \ref{sec:poset}.

\begin{theorem}\label{thm:invents} Let $T$ be an invertible tree and let $M$ be the unique perfect matching of $T$. If $A := A(X)$, then
\[
(A^{-1})_{a,b} = \begin{cases} 
(-1)^{1-m}, &\text{ if } a\neq b, \, |P_{a,b}| = 2m,\text{ and }P_{a,b} \text{ is $M$-alternating}; \\
0, &\text{ otherwise}. 
\end{cases}
\] 
\end{theorem}

\proof Since $X$ has an even number of vertices, $X\setminus a$ for any vertex $a$ has no perfect matching. We obtain by substituting $t=0$ into (\ref{eq:aa}) that $(-A)^{-1}_{a.a} = 0$ for all $a \in V(X)$. 

Let $a,b$ be distinct vertices of $X$. Evaluating (\ref{eq:ab}) at $t=0$ gives
\[
(-A)^{-1}_{a.b} = \frac{\phi(X\setminus P_{a,b} , 0)}{\phi(X,0)} .
\]
Since $X$ has a perfect matching, we have that $\phi(X,0) = (-1)^n$. By Lemma \ref{lem:alt}, $\phi(X \setminus P_{a,b}) \neq 0$ if and only if $P_{a,b}$ is a $M$-alternating path. If $P_{a,b}$ is a $M$-alternating path with $2m$ vertices, then $X \setminus P_{a,b}$ has a perfect matching and $\phi(X \setminus P_{a,b}) = (-1)^{n-m}$, 
whence the result follows. \qed 

Let $T$ be an invertible tree on $2n$ vertices with $A := A(T)$ and $M$ be the unique perfect matching of $T$. Let $(R,S)$ denote the bipartition of $T$. The above theorem implies that $A^{-1}$ has entries in $0, \pm 1$ and can be considered as the adjacency matrix of a signed graph $G^{\pm}$. Let $G$ denote the underlying graph of $G^{\pm}$. Observe that if $a,b$ are in the same part of the bipartition, $P_{a,b}$ has even number of edges and so $A^{-1}_{a.b} = 0$. Thus, $G$ is a bipartite graph and $(R,S)$ is also a bipartition of $G$. 

Let the edges of $M$ be $\{e_1, \ldots, e_n\}$. We may label the vertices of $T$ (and thus  $G$ and $G^{\pm}$ as they are defined on the same vertex set) as 
\[\{v_1,\ldots, v_n, w_1, \ldots, w_n\}\]
where $v_i$ is the end of $e_i$ in $R$ and $w_i$ is the end of $e_i$ in $S$. Let $\phi$ be the permutation of $V(T)$ which takes $v_i$ to $w_i$. Observe that $\phi(T)$ is an isomorphism of $T$ which fixes $M$. We will retain these definition for the rest of this section. 

\begin{lemma}\label{lem:edgeofinv} For $u,v \in V(T)$, $\phi(u)\phi(v)$ is an edge in $G$ if and only if $P_{u,v}$ is a $M$-alternating path in $T$.\end{lemma}

\proof This follows immediately from Theorem \ref{thm:invents}.\qed 

\begin{lemma}\label{lem:subgr} Let $T^{\pm}$ be the signed graph with $\phi(T)$ as its underlying graph, where the edges of $M$ are the positive edges and the edges not in $M$ are negative edges. Then $G^{\pm}$ contains $T^{\pm}$ as a signed subgraph; more precisely, $uv$ is an edge of $T^{\pm}$, then $uv$ is an edge of $G^{\pm}$ with the same sign.\end{lemma}

\proof Consider $v_i w_j$ an edge of $T^{\pm}$. We wish to show that $\phi(v_i w_j)$ is an edge of $G^{\pm}$ with the same sign as $v_i w_j$. If $v_i w_j \in M$, then we get that $(A^{-1})_{v_i w_j} = 1$ and $\phi(v_i w_j) = v_i w_j$ which has positive sign in $T^{\pm}$. Now suppose $v_i w_j \notin M$. We see that $\phi(v_i w_j) = \phi(v_i) \phi( w_j) = w_iv_j$. We have that
\[
P_{w_i,v_j} = \{w_i, e_i, v_i, v_iw_j, e_j, v_j\}
\]
is clearly a $M$-alternating path and so $(A^{-1})_{w_i v_j} = (-1)^{1-2} = -1$. \qed
 
Observe that $\phi(T)$ is a spanning tree of $G$. Each edge of $\phi(T)$ is one an unique edge-cut set of $G$ which contains no other edge of $\phi(T)$; these cuts are called the \textsl{fundamental cuts}. We say that a fundamental cut of $G$ (or of $G^{\pm}$) is \textsl{negative} if the corresponding edge of $\phi(T)$ is negative in $T^{\pm}$. 

\begin{lemma}\label{lem:fundcuts} Each edge $vw$ in $E(G^+) - E(\phi(T))$ is in $m-1$ negative fundamental cuts, where $2m-1$ is the distance from $v$ to $w$ in $\phi(T)$. \end{lemma}

\proof Let $Q$ be the unique path from $v$ to $w$ in $\phi(T)$. Consider $\phi(v)$ and $\phi(w)$ in $T$. Let $P_{\phi(v)\phi(w)}$ be the unique path in $T$ from $v$ to $w$.  Note that the proof of Lemma \ref{lem:subgr} actually implies that the negative edges of $\phi(T)$ are in one-to-one correspondence with the $M$-alternating paths in $T$ of length $3$. Since $vw$ is an edge of $G$ not in $\phi(T)$, we have that $P_{\phi(v)\phi(w)}$ is an $M$-alternating paths in $T$ of odd length, say $2m-1 \geq 5$.  The path $P_{\phi(v)\phi(w)}$ has the following sequence of vertices:
\[  \{u_0= \phi(v), u_1 = v, u_2, \ldots, u_{2m-3}, u_{2m-2} = w, u_{2m-1}= \phi(w)\}.
\]
Observe that $\{u_{2k}, u_{2k+1}, u_{2k+2}, u_{2k+3}\}$ for $k=0, \ldots, m-2$ are each  $M$-alternating paths in $T$ of length $3$ and $\phi(u_{2k}) \phi(u_{2k+3})$ is an edge of $T^{\pm}$ for each $k$. We have that $\phi(u_{2k}) \phi(u_{2k+3})$ for $k=0,\ldots, m-2$ and $u_{2k}u_{2k+1} \in M$ for $k= 1,\ldots, m-3$ are the edges of a path from $v$ to $w$ in $\phi(T)$. Since such a path is unique, these are the edges of $Q$. Then $Q$ has $m-1 + m-2 = m-3$ edges. The  $m-1$ edges of $Q$ of which are not in $\phi(T)$ and exactly the negative edges of $T'$ for which $vw$ is in the fundamental cut. \qed

For the following lemma, we note that we may consider any graph as a signed graph where the sign of every edge is positive. 

\begin{lemma}\label{lem:sgn} As signed graphs, $G^{\pm}$ and $G$ are switching equivalent. \end{lemma}

\proof Let $H$ be obtained from $G^{\pm}$ by switching on every negative cut of $G$. We claim that all edges of $H$ are positive, in which case $H$ would be equal to $G$ considered as a signed graph. 

It is clear that, in $H$, all edges of $\phi(T)$ now have positive sign. Consider $uv$ an edge of $G$ not in $\phi(T)$. The sign of $uv$ is $(-1)^{1-m}$ where $2m-1$ is the distance from $u$ to $v$ in $T$. By Lemma \ref{lem:fundcuts}, every negative edge is in an odd number of negative fundmental cuts and every positive edge is in an even number of negative fundmental cuts. After performing the switching described, every edge not in $\phi(T)$ will also have positive sign. \qed  

We have now reconstructed the theorem of Godsil. Note that in the original paper, the statement is imprecise about how $T$ is contained as a subgraph of its inverse graph. 

\begin{theorem}\label{thm:tree-inverse}\cite{Go85} If $T$ is a tree with a perfect matching $M$, then $A(T)^{-1}$ is a $(0,\pm 1)$ matrix, which can be considered as the adjacency matrix of a signed graph $G^{\pm}$. Further, $G^{\pm}$ is switching equivalent to a signed graph $G$ with all positive edges, which contains $\phi(T)$ as a spanning tree, where $\phi(T)$ is obtained from $T$ by the isomorphism which swaps the ends of every edge of $M$.   \end{theorem}

\proof This follows from Lemma \ref{lem:sgn}. \qed 

\begin{figure}[htbp]
\centering
\includegraphics[scale=0.67]{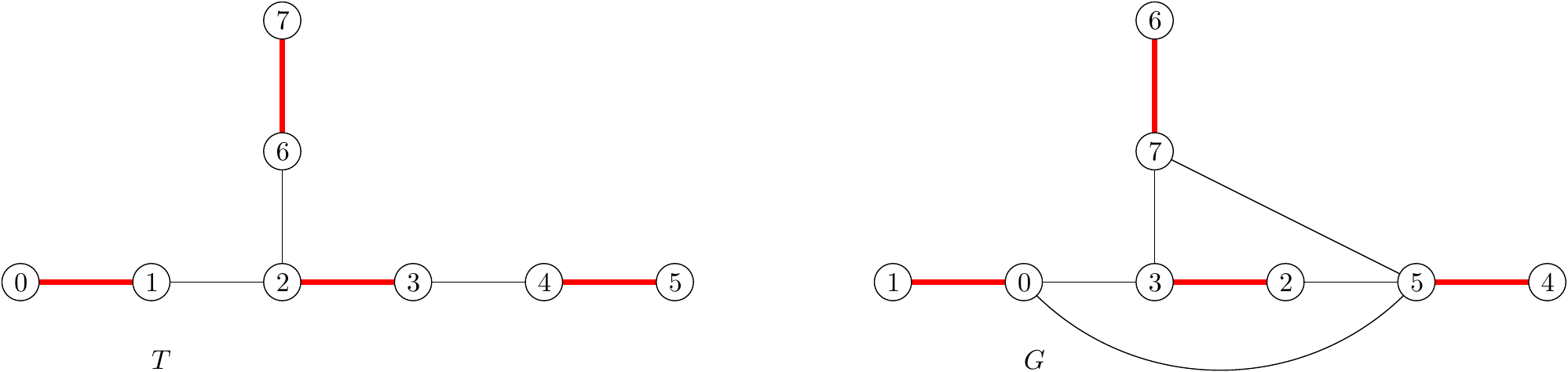}
\caption{An invertible tree $T$ and its inverse graph $G$ with the edges of the unique perfect matching of $T$ drawn as thicker edges in a lighter colour. \label{fig:ex}}
\end{figure}

Given an invertible tree $T$, we will refer to the graph $G$ of Theorem \ref{thm:tree-inverse}, forgetting the signature, as the inverse graph of $T$ and denote it $T^{-1}$. Figure \ref{fig:ex} shows an invertible tree on eight vertices and its inverse graph.

%
%
%

\section{A relation on invertible trees}\label{sec:poset}

Let $T$ be a tree with a unique perfect matching $M$. As in Section \ref{sec:inv-tree}, we will consider the inverse graph of $T$, the underlying graph of the signed graph given by $A(T)^{-1}$. We will denote the inverse graph of $T$ as $T^{-1}$. Recall the mapping $\phi$ is the involution permutation of the vertices of $T$ which swaps the ends of every edge of $M$. An edge of an invertible tree which is in the unique matching are said to be a \textsl{matching edge}, otherwise, it is a \textsl{non-matching edge}.

We will now define an operation on the class of invertible trees which maps a tree $T$ whose inverse graph $G$ has strictly more edges than $T$ to a tree whose inverse graph has fewer edges that $G$. Suppose $T$ is an invertible tree such that $T^{-1}$ contains an edge $e \notin \phi(T)$ and let $f$ be any non-matching edge of the fundamental cycle of $\phi(e)$ in $T \cup \{\phi(e)\}$. Let $\tau(T,e,f)$ be the tree obtained from $T$ by adding $\phi(e)$ and removing $f$. It is clear that the edges of $M$ remain a perfect matching in $X(T,e,f)$. If $\widetilde{T} = \tau(T,e,f)$ for some choice of edges $e$ and $f$, then we will say that $\widetilde{T}$ is obtained from $T$ by \textsl{tree-exchange}. 

\begin{figure}[htbp]
\centering
\includegraphics[scale=0.67]{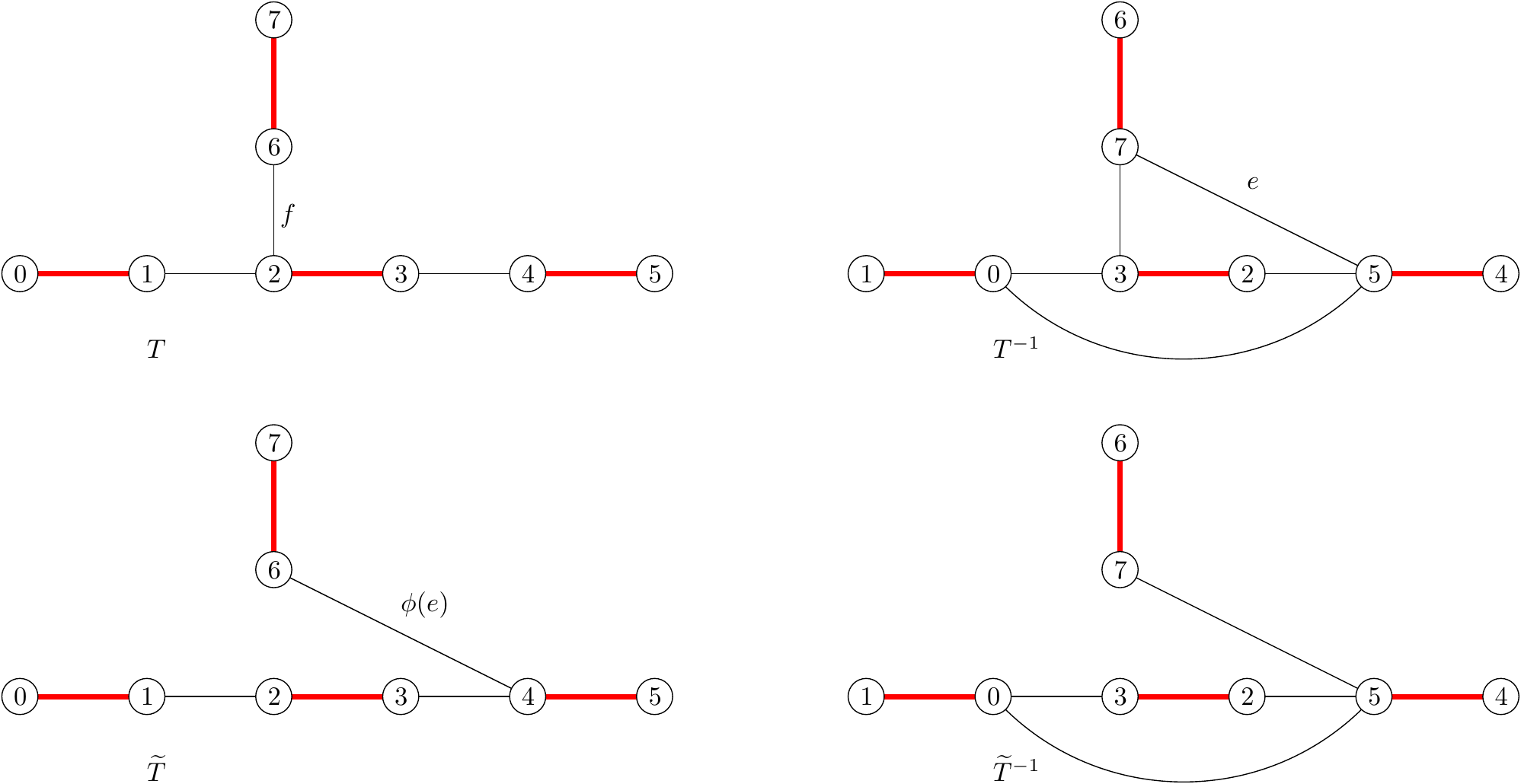}
\caption{In $T$ is an invertible graph and its inverse graph $T^{-1}$ with the edges of the unique perfect matching of $T$ drawn as thicker edges in a lighter colour. The tree $\widetilde{T}$ is obtained from $T$ by tree-exchange at edge $e$.\label{fig:ex2}}
\end{figure}

Figure \ref{fig:ex2} shows an example of tree-exchange; $T$ is an invertible tree whose inverse graph $T^{-1}$ contains $e = \{7,5\}$ such that $\phi(e) = \{6,4\}$ is not an edge of $T$. We have that $\widetilde{T} = \tau(T,e,f)$ where $f = \{6,2\}$. Note that $\phi(f) = \{7,3\}$ is not an edge of $\widetilde{T}^{-1}$.

\begin{lemma}\label{lem:switch-edges} If $T^{-1}$ contains an edge $e \notin \phi(T)$ and  $\widetilde{T}$ is obtained from $T$ by tree-exchange, then $\widetilde{T}^{-1}$ is a proper subgraph of $T^{-1}$.  \end{lemma}

\proof First we will show that $\widetilde{T}^{-1}$ is a subgraph of $T^{-1}$. We will then show that $\phi(f)$ is not an edge of $\widetilde{T}^{-1}$. Observe that $M$ is also the unique perfect matching of $\widetilde{T}$. 

Let $e = vw$. Since $e$ is an edge of $T^{-1}$, the path $P_{v,w}$ in $T$ is an $M$-alternating path. Since $e$ is not an edge of $T$, $P$ has order at least $5$. Observe that the matching edge at $v$ is $(v,\phi(v))$ and the matching edge at $w$ is $(w,\phi(w))$. Thus the first edge of $P_{v,w}$ is $(v,\phi(v))$  and the last edge of $P_{v,w}$ is $(w,\phi(w))$. We have that $C$, the fundamental cycle of $\phi(e)$ in $T \cup \{e\}$,  is the cycle consisting of $P_{v,w}$ with the first and last edge (that is   $(v,\phi(v))$ and $(w,\phi(w))$) deleted, together with $\phi(e)$. Observe that $C$ is not a $M$-alternating cycle, since both the edges preceding and following $\phi(e)$ are also non-matching edges. 

Consider an edge $uv$ of $\widetilde{T}^{-1}$. We many assume $uv$ is distinct from $e$, since we know $e$ is an edge of both graphs. The path $P_{u,v}$ from $u$ to $v$ in $\widetilde{T}$ is a $M$-alternating path. If $\phi(e)$ is not in $P_{u,v}$, then $P_{u,v}$ is also a $M$-alternating path in $T$ and so $uv$ is an edge of $T^{-1}$. Suppose $\phi(e)$ is an edge in $P_{u,v}$. Since $\phi(e)$ is not a matching edge, the edges immediately preceding and following $\phi(e)$ in $P_{u,v}$ are matching edges and must thus be $v\phi(v)$ and $w\phi(w)$. Let $Q_{\phi(v),\phi(w)}$ be the path on $C$ from $\phi(v)$ to $\phi(w)$  which does not use $\phi(e)$. If we replace $\phi(e)$ in $P_{u,v}$ with $Q_{\phi(v), \phi(w)}$, we obtain a walk $W$ in $T$ from $u$ to $v$. Observe that $Q$ is alternating, since the first and last edges of $Q_{\phi(v), \phi(w)}$ are non-matching edges. Since $T$ has no cycles, if $Q$ is not a walk, $Q$ must contain $\{a, ab, b, ba, a\}$ as a subsequence, for some vertices $a$ and $b$. This is not possible since $Q$ is $M$-alternating and thus there is a $M$-alternating path from $u$ to $v$ in $T$ and $uv$ is an edge of $T^{-1}$. 

Let $f = xy$. Now we will show that $\phi(x)$ is not adjacent to $\phi(y)$ in  $\widetilde{T}^{-1}$. Since $f$ is an edge of $T$, we have that $\phi(x)\phi(y) = \phi(f)$ is and edge of $T^{-1}$ and this would conclude the argument. Since $f$ is a non-matching edge on $C$, we have that $x\phi(x)$ and $y \phi(y)$ are also edges of $C$. In particular, we may take the path $P$ from $\phi(x)$ to $\phi(y)$ on $C$ which uses $\phi(e)$ but not $f$. Observe that $P$ does not contain $x\phi(x)$ or $y\phi(y)$ and is hence not a $M$-alternating path in $\widetilde{T}$. Since it is the unique path from $\phi(x)$ to $\phi(y)$, we obtain that $\phi(x)$ is not adjacent to $\phi(y)$ in $\widetilde{T}^{-1}$. \qed 

We obtain the following corollary.

\begin{corollary}\label{cor:lamddamed} If  $\widetilde{T}$ is obtained from $T$ by tree-exchange, then $\lambda_n(T) < \lambda_n(\widetilde{T})$.  \end{corollary}

\proof Since $\widetilde{T}^{-1}$ is a proper subgraph of  $T^{-1}$, we obtain that 
$\lambda_1( \widetilde{T}^{-1})$ is strictly less than $\lambda_1(T^{-1})$
which gives that 
\[
\frac{1}{\lambda_n(\widetilde{T})} < \frac{1}{\lambda_n(T)} 
\Leftrightarrow \lambda_n(T) < \lambda_n(\widetilde{T})
\]
where $2n$ is the number of vertices of $T$ and $\lambda_i(X)$ denotes the $i$th eigenvalue of $X$. \qed 

Let $\cT_{n}$ denote the class of invertible trees on $2n$ vertices. We define a relation $\leq_t$ on $\cT_n$ as follows:
\begin{enumerate}[(a)]
\item $T \leq_t T$ for all $T \in \cT_n$; and 
\item $T \leq_t \widetilde{T}$ if $\widetilde{T}$ is isomorphic to a tree obtained from $T$ by tree-exchange,
\end{enumerate}
and we extend by taking the transitive closure.

\begin{theorem} The relation $\leq_t$ is a partial order relation on $\cT_n$. \end{theorem}
\proof
Observe that if $T_1 \leq_t T_2$ and $T_1$ is not isomorphic to $T_2$, then by iterated applications of Corollary \ref{cor:lamddamed}, we see that $\lambda_n(T_1) < \lambda_n(T_2)$. This implies that $\leq_t$ is antisymmetric. Since $\leq_t$ is transitive by construction, it is a partial order on $\cT_n$. \qed
 
\begin{figure}[htbp]
\centering
\includegraphics[scale=0.85]{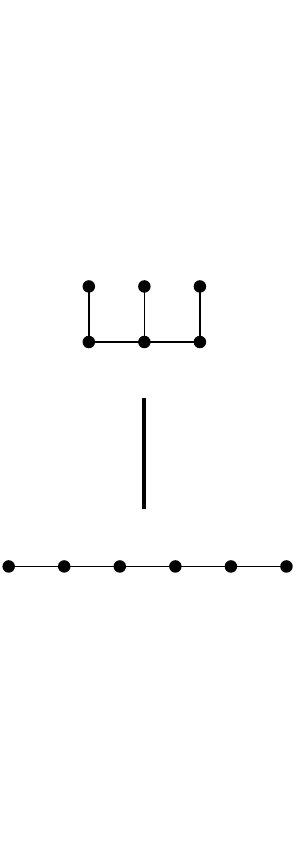}
\hspace{40pt}
\includegraphics[scale=0.8]{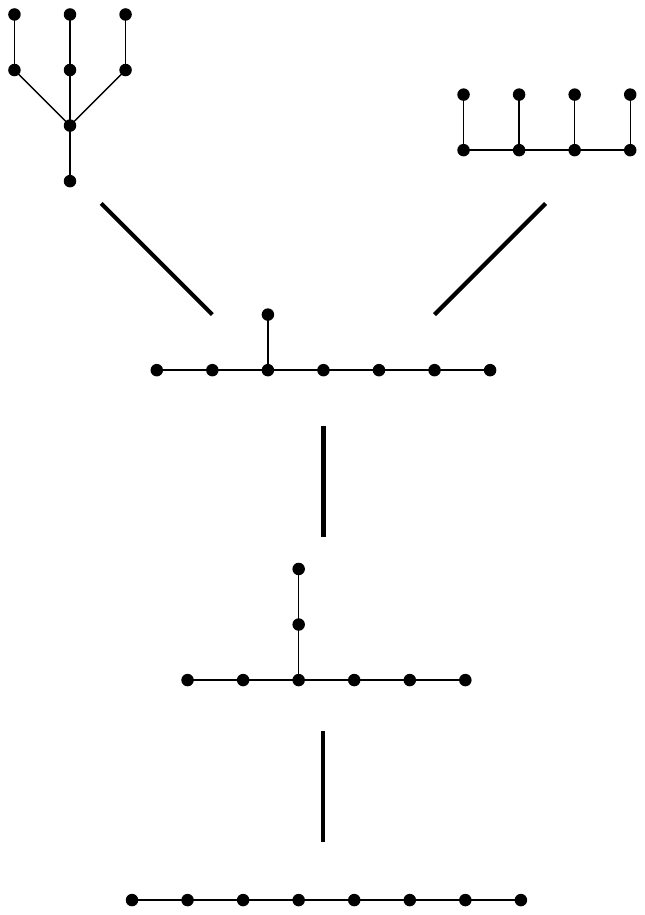}  
\caption{The Hasse diagrams of the invertible trees on $2n$ vertices, partially ordered by $\leq_t$, for $n = 3,4$.\label{fig:posets68}}
\end{figure}

We may consider $\cT_n$ as a poset under the partial ordering relation $\leq_t$. For $n \in \{1,2\}$, the only invertible trees on $2n$ vertices are the path graph. Figure \ref{fig:posets68} shows the Hasse diagrams of $\cT_3$ and $\cT_4$. 

\section{Maximal elements of $\cT_n$ ordered by $\leq_t$}

The maximal elements of $\cT_n$ ordered by $\leq_t$ are exactly the trees on $2n$ vertices which are isomorphic to their inverse graphs, which we will classify here. Note that this classification can be obtained using the result of \cite{GoMc78} on ``symmetric" characteristic polynomials of trees, but we will give direct proof here. 

Let $X$ be a graph with vertices $\{v_1,\ldots, v_n\}$ and let $Y$ be a disjoint union of rooted graphs $Y_1, \ldots, Y_n$. The \textsl{rooted product} of $X$ and $Y$, denoted $X(Y)$, is the graph obtain by identifying $v_i$ with the root vertex of $Y_i$. Figure \ref{fig:posets68} 

\begin{lemma}\label{lem:rtprod} If a tree $T$ on $2n$ vertices is isomorphic to its inverse graph, then $T$ is the rooted product of a tree on $n$ vertices and $n$ copies of rooted $K_2$.  \end{lemma} 

\proof We will show the statement by proving that every matching edge of $T$ has a vertex of degree $1$. Then, if all the vertices of degree $1$ are deleted, we obtain a tree on $n$ vertices, which shows that $T$ is the rooted product as described, where the matched edges form the $n$ copies of $K_2$, rooted at the vertex which does not have degree $1$ (which must exist, since $T$ is connected). 
 
 Let $A = A(T)$. We know that $A^{-1}$ is the adjacency matrix of the signing of $\phi(T)$ such that its unique matching, $M$, is exactly the set of positive edges. Let $P$ be the permutation matrix of $\phi$; observe that $P$ is also the adjacency matrix of the subgraph of $T$ with edge set $M$. We have that 
\[
A^{-1} = -A(\phi(T)) + 2P = - PAP + 2P = P(-AP + 2I).
\]
Thus $A^{-1}A = I$ gives  
\begin{equation}\label{eq:aainv} (-AP + 2I)A = P. \end{equation}
Let $N(u)$ denote the set of vertices adjacent to vertex $u$. Let $W(u)$ be the set of vertices matched to neighbours of $u$; that is to say,
\[
W(u) = \{\phi(w) \mid w \in N(u)\}.
\]
Note that $u \in W(u)$ for any vertex $u$. Noting that any vertex in $W(u)$ is matched to exactly one neighbour of $u$, we see that 
\[
(AP)_{u,v} = \begin{cases} 1, &\text{ if }v \in W(u); \\
0, &\text{ otherwise.}\end{cases}
\]
We obtain that
\[
(-AP + 2I)_{u,v}
= -1 |\{v\} \cap W(u)| + 2\delta_{u,v}
\]
where $\delta_{u,v}$ is the Kronecker delta.
Thus
\[
\begin{split}
((-AP + 2I)A)_{u,v} &= \sum_{y \in N(v)}  (-AP + 2I)_{u,y} \\
&= \sum_{y \in N(v)}  (-1) |\{y\} \cap W(u)| + 2\delta_{u,y} \\
&= (-1)|N(v) \cap W(u)| + 2\sum_{y \in N(v)} \delta_{u,y} \\
&= (-1)|N(v) \cap W(u)| + 2A_{u,v}. \\
\end{split}
\]
From (\ref{eq:aainv}) we obtain that 
\begin{equation}\label{eq:nvwu}
(-1)|N(v) \cap W(u)| + 2A_{u,v} = P_{u,v} = \delta_{u,\phi(v)}.
\end{equation}

Suppose there is a matched edge $w\phi(w)$ which is not incident to a vertex of degree $1$. Thus $w$ has a neighbour $u$ which is distinct from $\phi(w)$ and $\phi(w)$ has a neighbour $v$ which is distinct from $w$. Since there are no cycles in $T$, we have that $u\neq v$ and $u$ is not adjacent to $v$. We see from (\ref{eq:nvwu}) that 
\[
|N(v) \cap W(u)| =0.
\]
But this is a contradiction, since $\phi(w) \in N(v) \cap W(u)$. Thus, every matched edge is incident to a vertex of degree $1$, which completes the proof. \qed 

We see that $T \in \cT_n$ attaining the maximum $\lambda_n$ must be a maximal element of $\cT_n$ ordered by $\leq_t$.  Let $T_{2n}$ be the rooted product of the path $P_n$ with $n$ copies of $K_2$, each rooted at vertex $0$. This graph is sometimes called an \textsl{elongated caterpillar}. We can rederive the following result of [find citation] by showing that the elongated caterpillar attains the maximum $\lambda_n$ amongst all rooted product of trees on $n$ with $n$ copies of $K_2$. 

\begin{lemma} If $T$ is a tree on $2n$ vertices such that 
\[
\lambda_n(T) = \max_{X \in \cT} \lambda_n(X),
\]
then $T$ is isomorphic to $T_{2n}$. 
\end{lemma}

\proof By Lemma \ref{lem:switch-edges} and \ref{lem:rtprod}, we obtain that $T$ is a rooted product of a tree $Y$ with $K_2$. Theorem 2.1 of \cite{GoMc78}, we obtain that
\[ \begin{split}
\phi(A(T, t) &= \det(-tA(Y) + (t^2 - 1)I) \\
&= t^n\det\left(\frac{t^2-1}{t} I - A(Y)\right) \\
&= t^n\phi\left(Y, \frac{t^2 -1}{t}\right). 
\end{split}
\]
The roots of $\phi(A(T,t))$ are $\lambda$ such that $\frac{\lambda^2 -1}{\lambda} = \theta$, where $\theta$ is an eigenvalue of $Y$. For each eigenvalue $\theta$ of $Y$, we obtain two eigenvalues from the two roots of $t^2 - \theta t -1 =0$, which are
\[
\theta{\pm} = \frac{\theta \pm \sqrt{\theta^2 + 4}}{2}.
\]
The largest eigenvalue $\lambda_1(T)$ is $\theta_{+}$ where $\theta= \lambda_1(Y)$. Since $T$ is isomorphic to its inverse graph, we bet that $\lambda_n(T) = \frac{1}{\lambda_1(T)}$. Thus, $\lambda_n(T)$ is maximized amongst all graph $T$ which are isomorphic to its inverse graph when $\lambda_1(Y)$ is minimized, for all trees $Y$ on $n$ vertices. The tree with the smallest $\lambda_1(Y)$ is the path on $n$ vertices and so $T$ must be isomorphic to $T_{2n}$. 
 \qed 

 The eigenvalues of the paths are $2 \cos\left(\frac{\pi j}{n + 1}\right)$ for $ j = 1,\ldots, n$, and so the eigenvalues of $T_{2n}$ are $\lambda_j^{\pm}$ for $j = 1\, \ldots, n$, where 
\[
\cos\left(\frac{\pi j}{n + 1}\right) \pm \sqrt{\left(\cos\left(\frac{\pi j}{n + 1}\right)\right)^2+ 1}.
\]
Observe that $\lambda_1(T_{2n}) \leq 1 + \sqrt{2}$ and so $\lambda_n(T_{2n}) \geq \frac{1}{1+\sqrt{2}}$. 

\section{Minimal elements of $\cT_n$ ordered by $\leq_t$}

Godsil proved that the path $P_{2n}$ attains the minimum median eigenvalue in $\cT_n$. Thus $P_{2n}$ must be a minimal element of $\cT_n$ under $\leq_t$. We will show that the paths are the only minimal elements of $\cT_n$ under $\leq_t$.

\begin{lemma} If $T$ is a minimal element of $\cT_n$ under $\leq_t$, then $T$ has no vertex with degree at least $3$. 
\end{lemma}

\begin{figure}[htbp]
\centering
\includegraphics{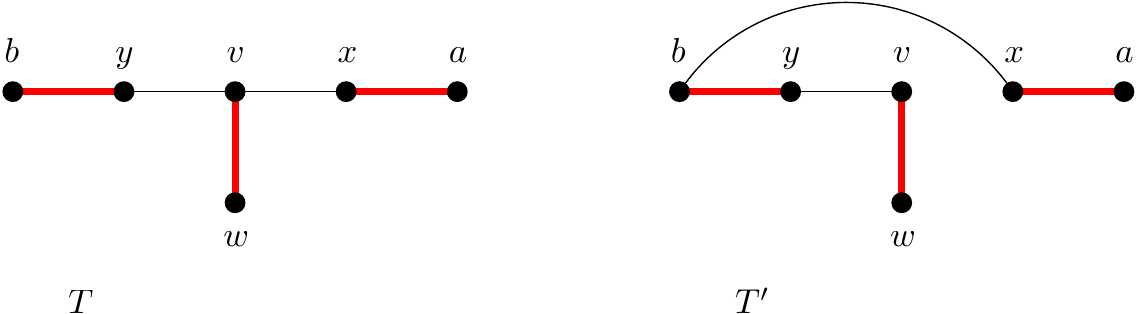}
\caption{the subgraphs of $T$ and $T'$ induced by $\{v,w,x,y,a,b\}$, with the edges of $M$ drawn as thicker edges in a lighter colour. \label{fig:deg3}}
\end{figure}

\proof Suppose for a contradiction that $T$ is a minimal element of $\cT_n$ under $\leq_t$ and $v$ is a vertex of $T$ with degree at least $3$. Let $M$ be the perfect matching of $T$ and $\phi$ be the mapping of the vertices of $T$ which switching the ends of the edges of $M$. Let $w$ be the neighbour of $v$ such that $vw \in M$. Let $x,y$ be two other neighbours of $v$ and let $a$ and $b$ be such that $xa, by \in M$. Let $T' = T \cup xb \setminus vx$. See Figure \ref{fig:deg3} for the subgraphs of $T$ and $T'$ induced by $\{v,w,x,y,a,b\}$. We will show that $T$ is obtained from $T'$ by tree-exchange. 

Observe that the path from $a$ to $w$ in $T'$ is a $M$-alternating path, and so $\phi(a)\phi(w)= xv$ is an edge of $T'^{-1}$ by Lemma \ref{lem:edgeofinv}. Thus $T = X(T',xv, xb)$, which is a contradiction. \qed 

\section{Open Problems} 

Tree-exchange allows us remove and add an edge from an invertible tree, such that the median eigenvalue strictly increases. One can ask if such an operation is possible for an bipartite graph with a unique matching. 

Considering the class of invertible trees as a poset, partially ordered by ``$\leq_t$'' as defined in Section \ref{sec:poset} gives rise to many questions. In particular, one may ask for a description of the covering relation of this poset and also which values  the Moebius function of this poset takes.
%


\end{document}